\documentclass[a4paper,10pt]{article}
\usepackage[utf8]{inputenc}
\usepackage{amsmath}
\usepackage{amsthm}
\usepackage[onehalfspacing]{setspace}
\usepackage{amssymb}
\usepackage{tikz}

\newtheorem{lem}{Lemma}

\newtheorem{thm}{Theorem}

\newtheorem{expl}{Example}

\title{A theorem about partitioning consecutive numbers}
\author{Kai Michael Renken}

\begin{document}

\maketitle

\begin{abstract}
In 1882 J.J. Sylvester already proved, that the number of different ways to partition a positive integer into consecutive positive integers exactly equals the number of odd divisors of that integer (see \cite{1}). We will now develop an interesting statement about triangular numbers, those positive integers which can be partitioned into consecutive numbers beginning at \(1\). For every partition of a triangular number \(n\) into consecutive numbers we can partition the sequence of numbers beginning at \(1\), adding up to \(n\) again, such that every part of this partition adds up to exactly one number of the chosen partition of \(n\).
\end{abstract}

\section{Preliminaries}
As an introducing example to the main theorem we want to prove later, consider the following "staircase shaped" Young tableau:
\begin{figure}[ht]
\centering
\begin{tikzpicture}[scale=0.5, line width=0.5pt]
  \draw (0,0) grid (5,1);
  \draw (0,0) grid (4,2);
  \draw (0,0) grid (3,3);
  \draw (0,0) grid (2,4);
  \draw (0,0) grid (1,5);
  
  \node[left] (A) at (1,0.5) {5};
  \node[left] (B) at (2,0.5) {5};
  \node[left] (C) at (3,0.5) {5};
  \node[left] (D) at (4,0.5) {5};
  \node[left] (E) at (5,0.5) {5};
  
  \node[left] (F) at (1,1.5) {4};
  \node[left] (G) at (2,1.5) {4};
  \node[left] (H) at (3,1.5) {4};
  \node[left] (I) at (4,1.5) {4};
  
  \node[left] (J) at (1,2.5) {3};
  \node[left] (K) at (2,2.5) {3};
  \node[left] (L) at (3,2.5) {3};
  
  \node[left] (M) at (1,3.5) {2};
  \node[left] (N) at (2,3.5) {2};
  
  \node[left] (O) at (1,4.5) {1};
  
\end{tikzpicture}
  \caption{"Staircase shaped" Young tableau consisting of 15 boxes}
  \label{figure1:Figure 1}
\end{figure}

We have one box in the first line and one more in every following line. Now, we can "rebuild" this tableau by reordering the lines, such that we have not necessarily one box in the first line anymore, but always have still one more box in every following line and such that we do not have to split the single lines, as the numbers show:

\begin{figure}[ht]
\centering
\begin{tikzpicture}[scale=0.5, line width=0.5pt]
  \draw (0,0) grid (8,1);
  \draw (0,0) grid (7,2);
  
  \node[left] (A) at (1,0.5) {5};
  \node[left] (B) at (2,0.5) {5};
  \node[left] (C) at (3,0.5) {5};
  \node[left] (D) at (4,0.5) {5};
  \node[left] (E) at (5,0.5) {5};
  \node[left] (F) at (6,0.5) {3};
  \node[left] (G) at (7,0.5) {3};
  \node[left] (H) at (8,0.5) {3};
  
  \node[left] (I) at (1,1.5) {4};
  \node[left] (J) at (2,1.5) {4};
  \node[left] (K) at (3,1.5) {4};
  \node[left] (L) at (4,1.5) {4};
  \node[left] (M) at (5,1.5) {2};
  \node[left] (N) at (6,1.5) {2};
  \node[left] (O) at (7,1.5) {1};

\end{tikzpicture}
  \caption{Rebuilt Young tableau}
  \label{figure2:Figure 2}
\end{figure}

The interesting question which will be answered by our main theorem is, if there is always a possibility to reorder the lines of a Young tableau of the first type without having to split them to build a Young tableau of the second type, provided that the numbers of boxes are the same. The question, how many different Young tableaux of this type exist, meaning Young tableaux consisting of a fixed number of boxes such that based on any line the following line has one box more, has already been answered in \cite{1}. Namely, for every odd divisor of the number of boxes, there exists exactly one Young tableau of this type.\\
Note, that in this paper we set \(\mathbb{N}\) to be the set of natural numbers without \(0\) and a sequence of consecutive numbers always denotes a sequence\\
\(n_1,\ldots,n_t\in\mathbb{N}\text{ (or }\mathbb{Z})\), such that \(n_{i+1}=n_i+1\). Furthermore, we use the notation \([n]:=1,2,3,\ldots,n\).
\begin{lem}\label{lemma1}
Let \(n,a,b\in\mathbb{N}\), such that \(n<a\leq b\) and \(\sum\limits_{i=1}^ni=\sum\limits_{i=a}^bi\). Then we have \(n\geq 2(b-a+1)\).
\begin{proof}
Note, that \(b-a+1\) is the number of summands in \(\sum\limits_{i=a}^bi\). Then we have:
\begin{align}
&\sum\limits_{i=1}^ni=\sum\limits_{i=a}^bi\notag\\
\Longleftrightarrow\quad &\frac{n(n+1)}{2}=(a+b)\frac{b-a+1}{2}\notag\\
\Longleftrightarrow\quad &\frac{n}{2}(n+1)=(b-a+1)\frac{a+b}{2}\notag
\end{align}
Now, by assumption we have \(\frac{a+b}{2}\geq n+1\), so we get \(\frac{n}{2}\geq b-a+1\), which is equivalent to \(n\geq 2(b-a+1)\).
\end{proof}
\end{lem}

\begin{lem}\label{lemma2}
For every \(m\in\mathbb{N}\) and every \(l\in\mathbb{Z}\) there exist pairs of numbers\\
\((x_1,x_1'),\ldots,(x_m,x_m')\in\mathbb{N}^2\), such that \(x_1,x_1',\ldots,x_m,x_m'\) are all distinct,\\
\(x_i'-x_i=i\) for every \(i=1\ldots,m\), \(l=\min\left\{x_1,x_1',\ldots,x_m,x_m'\right\}\) and\\
\(\max\left\{x_1,x_1',\ldots,x_m,x_m'\right\}\leq 2m+l\).
\begin{proof}
Consider the following pairs \((y_i,y_i')\) for \(i\) odd:
\[
\left(y_{2\left\lceil\frac{m}{2}\right\rceil-1},y_{2\left\lceil\frac{m}{2}\right\rceil-1}'\right)=\left(1,2\left\lceil\frac{m}{2}\right\rceil\right),\ldots,\left(y_1,y_1'\right)=\left(\left\lceil\frac{m}{2}\right\rceil,\left\lceil\frac{m}{2}\right\rceil+1\right)
\]
and the following for \(i\) even:
\begin{footnotesize}
\[
\left(y_{2\left\lfloor\frac{m}{2}\right\rfloor},y_{2\left\lfloor\frac{m}{2}\right\rfloor}'\right)=\left(2\left\lceil\frac{m}{2}\right\rceil+1,2m+1\right),\ldots,\left(y_2,y_2'\right)=\left(\left\lceil\frac{m}{2}\right\rceil+m,2m+2-\left\lfloor\frac{m}{2}\right\rfloor\right),
\]
\end{footnotesize}
where \((y_{i+2},y_{i+2}')=(y_i-1,y_i'+1)\) for all \(i=1,\ldots,m-2\). Obviously, these pairs satisfy the assumption \(y_i'-y_i=i\) and we have \(\min\left\{y_1,y_1',\ldots,y_m,y_m'\right\}=1\) and \(\max\left\{y_1,y_1',\ldots,y_m,y_m'\right\}=2m+1\) for every \(m\in\mathbb{N}\). (In fact, we even have \(\max\left\{y_1,y_1'\right\}=2m<2m+1\) for the case \(m=1\)). Now, set \(x_i:=y_i+l-1\) and \(x_i'=y_i'+l-1\) for all \(i=1,\ldots,m\) and we are done.
\end{proof}
\end{lem}

\section{The main theorem}

\begin{thm}\label{theorem1}
Let \(n,a,b\in\mathbb{N}\) (\(b\geq a\)), such that \(\sum\limits_{i=1}^ni=\sum\limits_{i=a}^bi\), then for every \(a\leq t\leq b\) there exists a subset \(U_t\subseteq [n]\), such that \(U_i\cap U_j=\emptyset\) for all \(i\neq j\),\\
\([n]=\bigcup\limits_{a\leq t\leq b}U_t\) and \(\sum\limits_{i\in U_t}i=t\).
\begin{proof}
Without loss of generality let \(a>n\). Otherwise, we have \(U_t=\{t\}\) for all \(a\leq t\leq n\) and \(\sum\limits_{i=1}^{a-1}i=\sum\limits_{i=n+1}^bi\), so the remaining problem is reduced to the requested case, since we have \(n+1>a-1\).\\
Now, denote the number of summands in the sum \(\sum\limits_{i=a}^bi\) by \(s:=b-a+1\) and set:
\begin{align}
P&=\left\{p_1,\ldots,p_s\right\}:=\left\{n-2s+1,n-2s+2,\ldots,n-s-1,n-s\right\},\notag\\
Q&=\left\{q_1,\ldots,q_s\right\}:=\left\{n-s+1,n-s+2,\ldots,n-1,n\right\},\notag\\
R&=\left\{r_1,\ldots,r_s\right\}:=\left\{a,a+1,\ldots,b-1,b\right\}\notag
\end{align}
We have \(P,Q,R\subset\mathbb{N}\), since \(n-2s+1>0\) holds by Lemma \ref{lemma1}. Furthermore, we see that we have \(p_i+q_j=2n-2s+1\) for all \(i,j\), such that \(i+j=s+1\). Now, set \(c:=2n-2s+1\) and for every \(t\) such that \(r_t-c<0\) consider pairs \(\left\{p_{i_t},q_{j_t}\right\}\) and \(\left\{p_{i_t'},q_{j_t'}\right\}\), such that \(i_t+j_t=i_t'+j_t'=s+1\), \(q_{j_t'}-q_{j_t}=c-r_t\) and all appearing numbers in all pairs are distinct (we actually only have to require that all \(q_{j_t},q_{j_t'}\) are distinct, then all the rest is distinct by the condition \(i_t+j_t=i_t'+j_t'=s+1\)). We can find those kinds of tuples of pairs for every \(t\) satisfying \(r_t-c<0\) by Lemma \ref{lemma2} as follows. For every \(t\), such that \(r_t-c<0\) there exists a \(t'\), such that \(r_{t'}-c=-(r_t-c)\), since otherwise we had \(\sum\limits_{i=1}^sp_i+q_i>\sum\limits_{i=1}^sr_i\), because the numbers \(r_1-c,\ldots,r_s-c\) are obviously consecutive, but this is a contradiction to the assumption \(\sum\limits_{i=1}^ni=\sum\limits_{i=a}^bi\). By this observation, there must still exist a \(t\), such that \(r_t-c=0\), since \(r_1,\ldots,r_s\) are consecutive numbers, so letting \(m\) denote the maximum number, such that there exists a \(t\) satisfying \(r_t-c=-m\), we have at least \(2m+1\) consecutive numbers in \(Q\). Now, choose such a set of \(2m+1\) consecutive numbers from \(Q\) and let \(l\) denote the minimum number in this set, then applying Lemma \ref{lemma2} we get the requested tuples of pairs satisfying \(q_{j_t'}-q_{j_t}=c-r_t\). Finally, set \(U_{r_t}:=\left\{p_{i_t'},q_{j_t}\right\}\) and \(U_{r_t+2(c-r_t)}:=\left\{p_{i_t},q_{j_t'}\right\}\) for every \(t\) satisfying \(r_t-c<0\) and we have \(\sum\limits_{i\in U_{r_t}}i=r_t\) and \(\sum\limits_{i\in U_{r_t+2(c-r_t)}}i=r_t+2(c-r_t)\). For the remaining \(t\)'s satisfying \(r_t-c\geq 0\) we set \(U_{r_t}:=\left\{p_{i_t},q_{j_t}\right\}\), where \(\left\{p_{i_t},q_{j_t}\right\}\) can be any of the remaining pairs, satisfying \(p_{i_t}+q_{j_t}=c\). For those \(t\)'s we have \(\sum\limits_{i\in U_{r_t}}i=c\). So, the problem is reduced to a smaller one of the type \(\sum\limits_{i=1}^{n-2s}i=\sum\limits_{i=r_k-c}^{r_s-c}i\), where \(k\) is the minimal number such that \(r_k-c>0\). By induction we are done as follows. The sets \(U_{r_t}\) and \(U_{r_t+2(c-r_t)}\), where \(r_t-c<0\) stay the same till the end, whereas the other sets \(U_{r_t}\), where we still have \(\sum\limits_{i\in U_{r_t}}i=c<r_t\) will be "filled up" during the next steps of induction.
\end{proof}
\end{thm}
The reader may develop a better intuition for the preceding proof by considering the following example.
\begin{expl}
Consider \(n=14\), \(a=15\) and \(b=20\), then we obviously have \(\sum\limits_{i=1}^ni=\sum\limits_{i=a}^bi\). Now, according to the proof of Theorem \ref{theorem1} we have \(s=6\) and we consider the following pairs:
\begin{align}
&\left\{p_1,q_6\right\}=\left\{3,14\right\},\quad \left\{p_2,q_5\right\}=\left\{4,13\right\},\quad \left\{p_3,q_4\right\}=\left\{5,12\right\},\notag\\
&\left\{p_4,q_3\right\}=\left\{6,11\right\},\quad \left\{p_5,q_2\right\}=\left\{7,10\right\},\quad \left\{p_6,q_1\right\}=\left\{8,9\right\}\notag
\end{align}
We see that we have \(c=17\), so by subtracting \(c\) from the numbers \(15,\ldots,20\) only the last 3 numbers stay strictly positive. Now we "swap" the second components of two pairs just as we did in the preceding proof and construct the following sets:
\begin{align}
&U_{15}=\left\{3,12\right\},\quad U_{16}=\left\{6,10\right\},\quad U_{17}=\left\{8,9\right\},\notag\\
&U_{18}=\left\{7,11\right\},\quad U_{19}=\left\{5,14\right\},\quad U_{20}=\left\{4,13\right\}\notag
\end{align}
We can see that we already have \(\sum\limits_{i\in U_t}i=t\) for all \(t=15,\ldots,19\) and\\
\(\sum\limits_{i\in U_{20}}i=c=17<20\), so the problem is reduced to \(\sum\limits_{i=1}^2i=3\) and inductively we will be done by the next step. Formally this means, that \(U_{15},\ldots,U_{19}\) stay the same and we "fill up" the remaining set \(U_{20}\) with the remaining numbers \(1,2\), such that we get \(U_{20}=\left\{13,4,1,2\right\}\).
\end{expl}
There are obviously many more ways of partitioning numbers in our sense than the one way the "algorithm" in the proof of Theorem \ref{theorem1} gives us, but until now this is the only working way we know in general and we even still do not know, how many partitions leading to the requested result exist.

\end{document}